\documentclass{article}
\usepackage{amsmath}
\usepackage{amssymb}
\usepackage{latexsym}
\usepackage{verbatim}
\usepackage{graphics}
\usepackage{amsfonts}
\usepackage{graphicx}
\usepackage{epstopdf}
\input xypic
\newtheorem{thm}{Theorem}
\newtheorem{conj}{Conjecture}
\newtheorem{prop}{Proposition}

\newtheorem{cor}{Corollary}

\begin{document}

\title{An Introduction to Virtual Spatial Graph Theory}

\author{Thomas Fleming\\
    Department of Mathematics\\
    University of California, San Diego\\
    La Jolla, CA 92093-0112\\
    {\it  tfleming@math.ucsd.edu}
   \and
    Blake Mellor\\
    Mathematics Department\\
    Loyola Marymount University\\
    Los Angeles, CA  90045-2659\\
    {\it  bmellor@lmu.edu}}

\date{}

\maketitle

\begin{abstract}

Two natural generalizations of knot theory are the study of spatial
graphs and virtual knots. Our goal is to unify these two approaches
into the study of virtual spatial graphs.  We state the definitions,
provide some examples, and survey the known results.  We hope that
this paper will help lead to rapid development of the area.

\end{abstract}



\section{Introduction} \label{S:intro}


The mathematical theory of knots studies the many ways a single loop
can be tangled up in space.  Since many biological molecules, such
as DNA, often form loops, knot theory has been applied to biological
systems with good effect \cite{dw}.  However, many biological
molecules form far more complicated shapes than simple loops;
proteins, for example, often contain extensive crosslinking between
cystine residues, and hence from the mathematical viewpoint are far
more complicated structures--spatial graphs.  The study of graphs
embedded in space is known as spatial graph theory, and researchers
such as Flapan \cite{fl} have obtained good results by applying it
to chemical problems. However, in biological systems, proteins are
often associated with membranes, meaning that some portions of the
molecule are prevented from interacting with others.  In the case of
a simple loop, the virtual knot theory of Kauffman \cite{ka}
provides a mathematical framework for studying such systems, as it
allows some crossings of strands to be labeled ``virtual,'' i.e.
non-interacting.  We hope that a merging of these two theories,
called virtual spatial graph theory, will prove equally useful in
the biological sciences.

Knot theory studies embeddings of circles up to isotopy. There are
many ways to extend the ideas of knot theory; two natural choices
are the study of spatial graphs and the theory of virtual knots. The
theory of spatial graphs generalizes the objects we embed, by
studying isotopy classes of arbitrary graphs embedded in $S^{3}$.
The theory of virtual knots is quite different--it generalizes the
idea of an embedding. Knots and links can be studied by projecting
the embedding into a plane, but retaining information about over-
and under-crossings. These knot diagrams are then taken up to an
equivalence defined by Reidemeister moves.  Virtual knot theory
simply allows a third type of crossing, a ``virtual'' crossing, and
introduces new Reidemeister moves which determine how this type of
crossing behaves.  Virtual knot theory then studies these virtual
knot diagrams up to equivalence under the full set of Reidemeister
moves.  These two generalizations of knot theory move in very
different directions.  It is our intention to find a unified
approach.

$$ \xymatrix@1{ &
{\rm Virtual~Spatial~Graph~Theory}
\\ \\
{\rm Spatial~Graph~Theory} \ar@{^{(}-->}[uur] & & {\rm
Virtual~Knot~Theory} \ar@{_{(}->}[uul]
\\
\\
 & {\rm Knot~Theory} \ar@{_{(}->}[uul] \ar@{^{(}->}[uur]
 }$$

Naturally, spatial graphs may be studied by examining their
diagrams, and by introducing virtual crossings and new Reidemeister
moves, we may define virtual spatial graphs. A \emph{virtual spatial
graph} is an immersion of a graph $G$ into the plane, with crossings
labeled over, under, or virtual.  These diagrams are taken up to the
equivalence relation generated by the Reidemeister moves of Figures
\ref{F:classicalmoves} and \ref{F:virtualmoves}.  Many invariants of
both virtual knots and spatial graphs may be extended to these new
objects.  This is the topic of Section \ref{S:invariants}.

A knot can be described by the Gauss code of its projection--the
sequence that records the order the crossings are met as we traverse
the knot diagram. However, there are many more such sequences than
there are real knots; those that correspond to classical knots are
known as ``realizable'' codes.  One motivation for virtual knots is
that they provide realizations for the Gauss codes that do not
correspond to classical graphs.  A similar motivation for the study
of virtual spatial graphs comes from the Gauss code of a diagram for
a spatial graph.  The Gauss code for a spatial graph simply records
the sequence of crossings along each edge of the graph.  Again, not
every such code corresponds to a classical spatial graph, but any
such code corresponds to some virtual spatial graph.  In Section
\ref{gcodes} we will define a Gauss code for spatial graphs and
determine which Gauss codes can be realized by a classical spatial
embedding.

Another motivation for virtual spatial graph theory is that it can
give insights into classical problems in topological graph theory.
In Section \ref{S:ivl} we will use virtual spatial graph theory to
produce a filtration on the set of intrinsically linked graphs.
Sachs, and Conway and Gordon proved that every embedding of $K_{6}$
into $S^{3}$ contains a non-split link \cite{sachs, c&g}.  A graph
with this property is called \emph{intrinsically linked}, and these
graphs have been extensively studied.  Given two such graphs, one
might want to compare their ``linkedness.''  One measure of this
property might be the minimum number of links contained in any
spatial embedding.  For example, there is an embedding of $K_{6}$
which contains exactly one non-split link, but every embedding of
$K_{4,4}$ contains at least two non-split links \cite{fm4}. However,
both $K_{6}$ and the Petersen graph have embeddings with exactly one
link.  Yet, in a certain sense, $K_{6}$ is more linked than the
Petersen graph. That is, we say a graph $G$ is
\emph{$n$-intrinsically virtually linked} ($n$-IVL) if every virtual
spatial graph diagram of $G$ with $n$ or fewer virtual crossings
contains a nonsplit virtual link. This produces a filtration of the
set of intrinsically linked graphs, and the deeper in the filtration
a graph persists, the more ``linked'' the graph is.  As we will
discuss in Section \ref{S:ivl}, $K_{6}$ is 2-IVL, but the Petersen
graph is not.

Virtual spatial graph theory is relatively new, and many aspects of the theory are in early stages
of development.  Thus, there are many accessible open problems and questions in this field, and
we hope it remains an active area of research for some time.

\smallskip

\noindent{\sc Acknowledgement:} The authors would like to thank Akio
Kawauchi and the organizers of the conference \emph{Knot Theory for
Scientific Objects} for giving them the opportunity to introduce
this work.

\section{Virtual Spatial Graphs} \label{S:definition}

We will use a combinatorial definition of virtual spatial graphs. A
graph is a set of vertices V and a set of edges $E
\subset V \times V$.  Unless otherwise stated, we will consider vertex-oriented and directed
graphs, so that each edge is an ordered pair of
vertices.  A spatial graph is an embedding of $G$ in $S^3$ that maps
the vertices to points and an edge $(u, v)$ to an arc whose
endpoints are the images of the vertices $u$ and $v$, and that is
oriented from $u$ to $v$.  We will consider these embeddings up to
ambient isotopy.  We can always represent such an embedding by
projecting it to a plane so that each vertex neighborhood is a
collection of rays with one end at the vertex and crossings of edges
of the graph are transverse double points in the interior of the
edges (as in the usual knot and link diagrams).

\begin{figure}
    $$\includegraphics{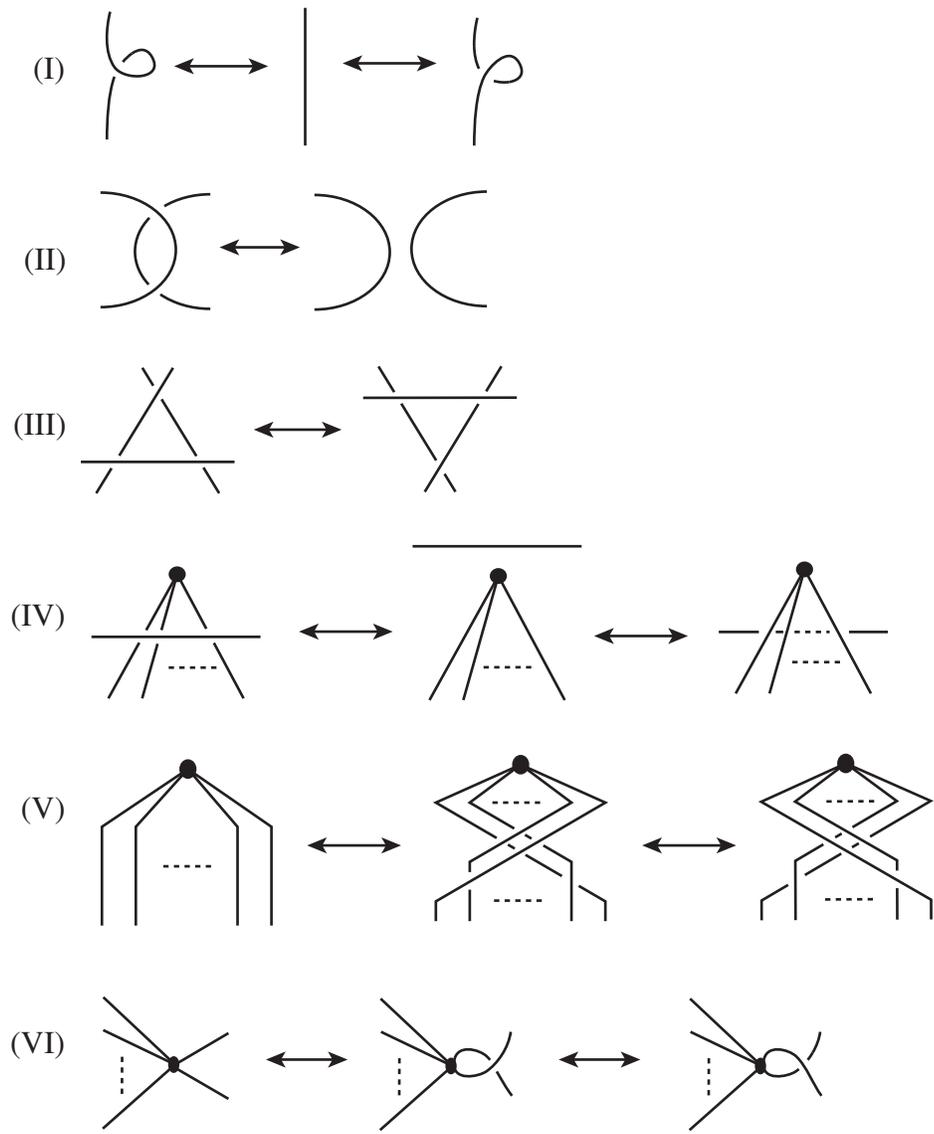}$$
    \caption{Reidemeister moves for graphs} \label{F:classicalmoves}
    \end{figure}

Kauffman \cite{ka2} and Yamada \cite{ya} have shown that ambient
isotopy of spatial graphs is generated by a set of local moves on
these diagrams which generalize the Reidemeister moves for knots and
links.  These Reidemeister moves for graphs are shown in
Figure~\ref{F:classicalmoves}.

\begin{figure}
    $$\includegraphics{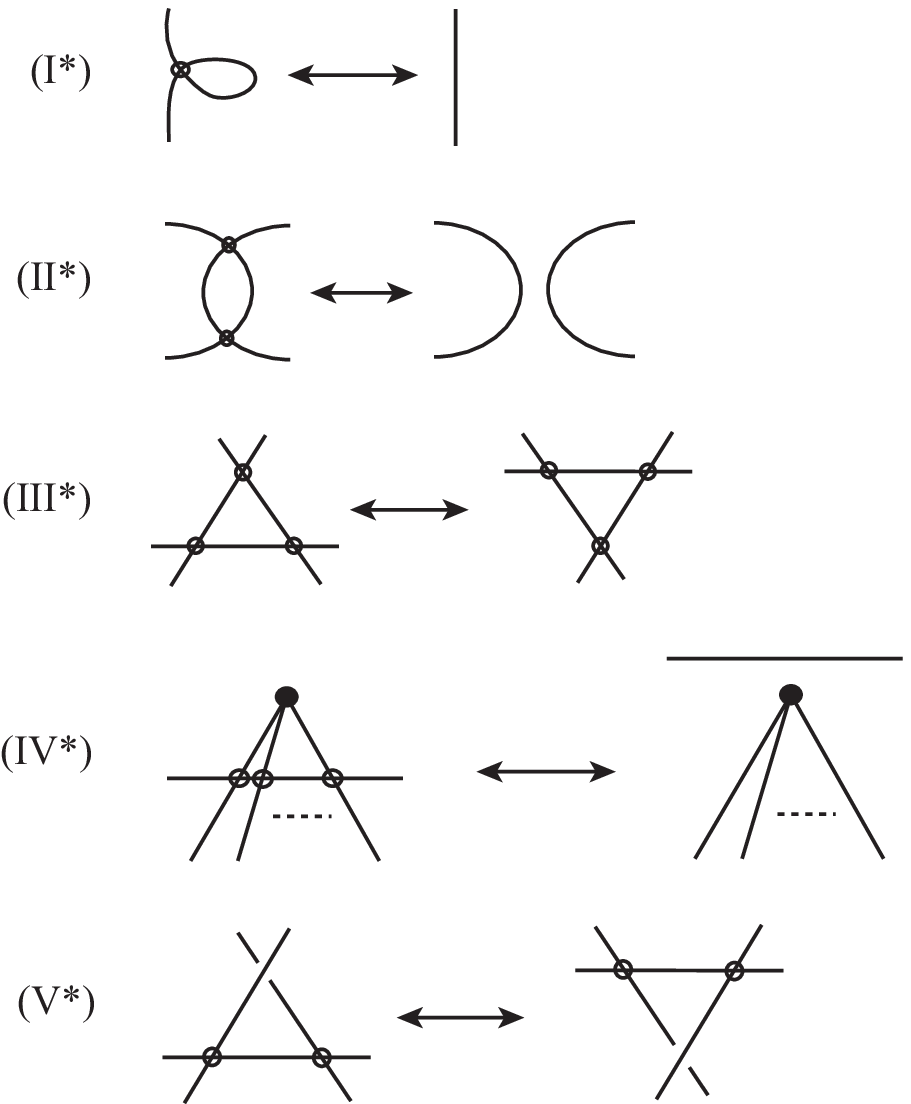}$$
    \caption{Reidemeister moves for virtual graphs} \label{F:virtualmoves}
    \end{figure}

A {\it virtual graph diagram} is a classical graph
diagram, with the addition of virtual crossings. The idea is that
the virtual crossings are not really there (hence the name
``virtual"). To make sense of this, we extend our set of
Reidemeister moves for graphs to include moves with virtual
crossings.  We need to introduce five more moves, (I*) - (V*), shown in
Figure~\ref{F:virtualmoves}.  Notice that moves (I*) - (IV*) are
just the purely virtual versions of moves (I) - (IV); move (V*) is
the only move which combines classical and virtual crossings (in
fact, there are two versions of the move, since the classical
crossing may be either positive or negative).

\begin{figure}
    $$\includegraphics{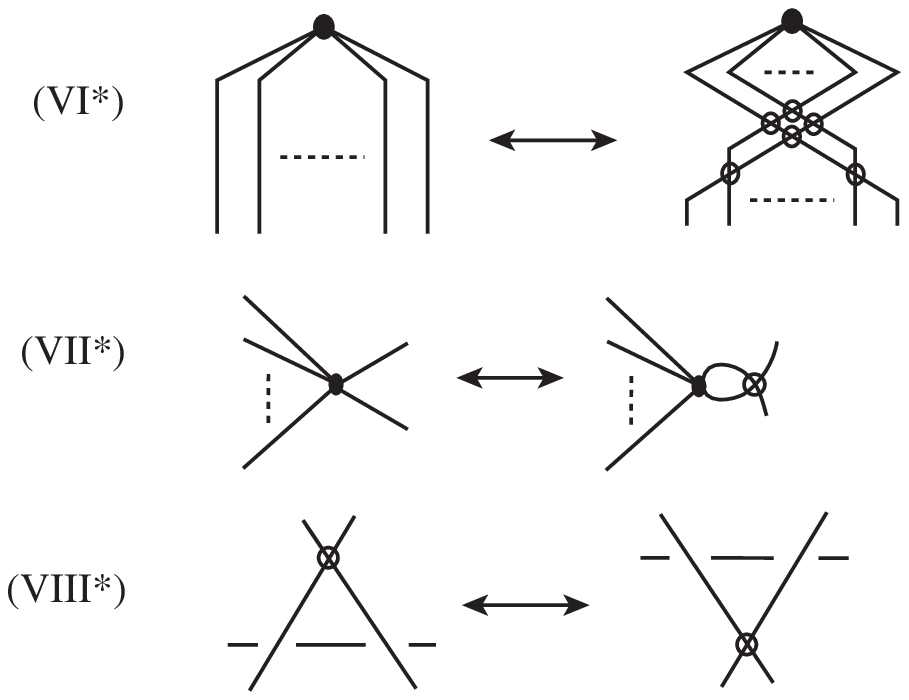}$$
    \caption{Forbidden Reidemeister moves for virtual graphs} \label{F:forbiddenmoves}
    \end{figure}

There are also three moves which, while they might seem reasonable, are {\it not} allowed.
These {\it forbidden moves} are shown in Figure~\ref{F:forbiddenmoves}.

\subsection{Forbidden Moves} \label{SS:forbidden}

If we allow the forbidden moves in Figure \ref{F:forbiddenmoves},
then many more virtual graph diagrams become equivalent.  In the
case of knots, allowing move (VIII*) trivializes the theory, and all
virtual knots become trivial \cite{ka, ne}.  However, when we look
at virtual links or virtual graph diagrams, the effect is not quite
so drastic.  The following proposition shows that the forbidden
moves do not trivialize virtual graph theory, as they do virtual
knot theory.

\begin{prop} \label{P:linkingnumber}
There are virtual graph diagrams (of a connected graph) which are not equivalent modulo the
forbidden moves.
\end{prop}


We use the virtual spatial graph invariant $T(G)$, defined in
Section \ref{S:invariants} below. The pairwise linking numbers for all the links
in $T(G)$ can be computed by using the Gauss formula
($\frac{1}{2}$(number of positive crossings) - $\frac{1}{2}$(number
of negative crossings)). If the links are virtual, these linking
numbers may be half-integers, but they are still invariant under all
the classical and virtual Reidemeister moves and the forbidden moves
(VI*), (VII*) and (VIII*).

The two virtual graph diagrams on the right in Figure \ref{F:Texample} have links in
$T(G)$ with different linking numbers, and so the diagrams are inequivalent, even
allowing the forbidden moves.


\subsection{Invariants of Virtual Spatial Graphs}
\label{S:invariants}

\begin{figure} [ht]
    $$\includegraphics{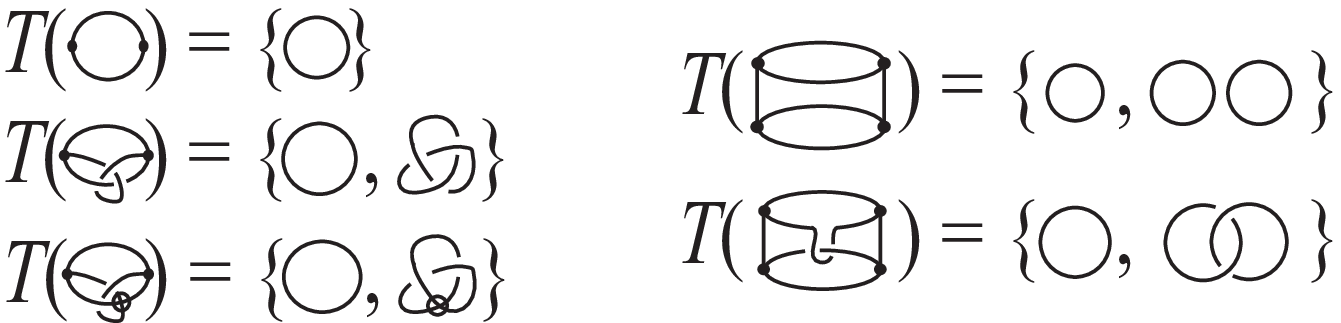}$$
    \caption{Examples of T(G) for virtual graph diagrams} \label{F:Texample}
    \end{figure}

Kauffman \cite{ka2} introduced a topological invariant of a spatial
graph defined as the collection of all knots and links formed by a
local replacement at each vertex of the graph.  Each local
replacement joins two of the edges incident to the vertex and leaves
the other edges as free ends (i.e. creates new vertices of degree
one at the end of each of the other edges). Choosing a replacement
at each vertex of a graph $G$ creates a link $L(G)$ (after erasing
all unknotted arcs).  $T(G)$ is the collection of all links $L(G)$
for all possible choices of replacements. For virtual graphs, we can
define $T(G)$ in exactly the same way, except that it is now a
collection of {\it virtual} links. Examples are shown in Figure \ref{F:Texample}.

The {\it fundamental group} of a classical knot or spatial graph is
the fundamental group of its complement in $S^3$.  Given a diagram
for the knot or graph, this group can be calculated using the
{\it Wirtinger presentation}. Kauffman \cite{ka} defined the fundamental group of a virtual knot
by modifying the Wirtinger presentation for classical knots.  We define the fundamental group of
a
virtual spatial graph in the same way, by modifying the Wirtinger presentation for classical
spatial graphs.  Note that for a classical spatial graph, this modified calculation
still produces the usual fundamental group.


The quandle is a combinatorial knot invariant introduced by Joyce \cite{joyce}, that
was generalized to virtual knots by Kauffman \cite{ka}, and
strengthed by Manturov \cite{ma}. Modifying Manturov's approach, we
can construct a similar invariant for virtual spatial graphs, though
for general graphs this invariant is less potent than in the case of
knots.  Information about the embedding
is encoded by relations generated by the crossings and vertices.



Yamada introduced a polynomial invariant $R$ of spatial graphs in
\cite{ya}. Using skein relations, $R(G)$ can be computed by reducing
the graph $G$ to a bouquet of circles, where a (classical) trivial
bouquet of $n$ circles has $Y(B(n)) = -(-\sigma)^{n}$. In the case
of virtual spatial graphs we can use exactly the same skein
relations to compute $R(G)$, simply by ignoring virtual crossings.
The only difference is that we may end up with a {\it virtual
bouquet}--a bouquet of circles with only virtual crossings. If $G$
is a virtual bouquet of $n$ circles, we simply define $R(G) = R(B_n)
= -(-\sigma)^n$. This gives a virtual graph diagram invariant.  It
is easy to find nontrivial virtual graph diagrams of planar graphs
with trivial Yamada polynomial.

Detailed discussion of these invariants can be found in
\cite{fm2}.


\section{Gauss Codes}
\label{gcodes}

\begin{figure} [ht]
    $$\includegraphics{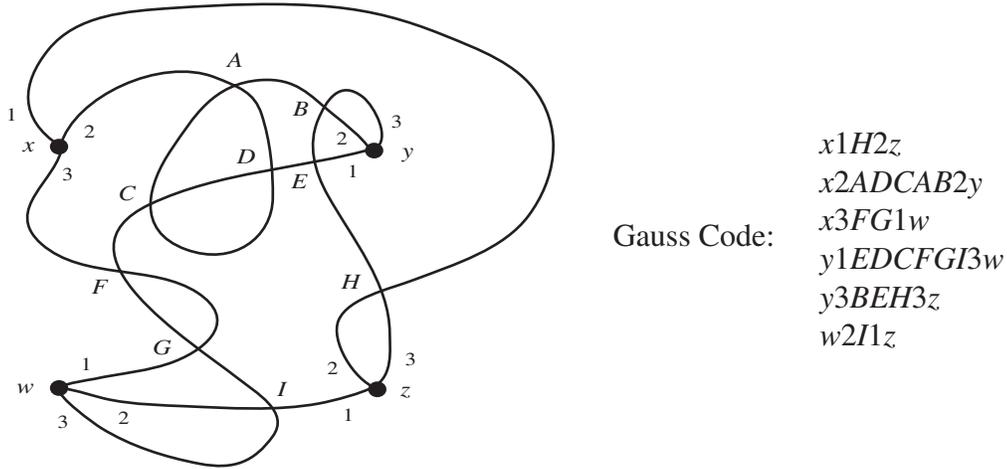}$$
    \caption{Gauss code for the shadow of a graph diagram} \label{F:gausscode}
    \end{figure}

It is easiest to associate a Gauss
code to an immersion of a closed curve (or graph) in the plane, so
we first study the shadow of our graph diagram, where the over/under
information at the crossings is ignored.
Figure~\ref{F:gausscode} illustrates the Gauss
code for such a shadow.
To produce the Gauss code for the original diagram, we augment the
Gauss code for its shadow by recording whether each crossing is an
over-crossing or an under-crossing.  If the graph is
directed, we can also label each crossing by its sign.

An important problem in the study of Gauss codes is to find
algorithms for determining whether a Gauss code is realizable, that
is, the Gauss code of a classical embedding.  For closed curves,
there are several algorithms \cite{fo, de, rr}; these methods can be
generalized to arbitrary graphs \cite{fm}.

Given a Gauss code $S$ for an abstract graph $G$, there is a split code $S^{*}$ produced in a
similar manner as when studying Gauss codes for a knot.  However, in this process one must also
alter the graph $G$ to produce a chord diagram on a \emph{split graph} $G^{*}$.  That is, given
a Gauss code $S$ on $G$, we may label each edge of $G$ with a sequence of symbols.  When
splitting $S$ at symbol $b$, we cut the edges of $G$ at the labels $b$ and then reconnect them to
produce a connected graph.  We also add a chord running from one instance of $B$ to the other.

\begin{thm}
\label{codethm}
Let $S$ be a Gauss code for an abstract connected
graph $G$, with split code $S^*$ and split graph $G^*$.  Then $S$ is
realizable by a classical spatial graph embedding if and only if
$G^*$ is planar, and the embedding of $G^*$ determined by $S^*$ can
be extended to an embedding of the corresponding chord diagram.
\end{thm}

Perhaps surprisingly, checking the conditions of Theorem
\ref{codethm} is computationally easy.

\begin{prop}
Given a Gauss Code $S$ on a connected graph $G$, there is a
polynomial time algorithm to determine whether $S$ is realizable by a classical embedding of
$G$.
\end{prop}

This algorithm and a proof of Theorem \ref{codethm} are discussed in
\cite{fm}.

Virtual graph diagrams also have Gauss codes, produced in exactly
the same way, except that virtual crossings are ignored.  As Theorem \ref{T:virtualrealize}
shows, one motivation for
studying virtual graph diagrams is that they allow us to realize the
``unrealizable" Gauss codes.

\begin{thm} \label{T:virtualrealize}
Every Gauss code can be realized as the code for a virtual graph
diagram.
\end{thm}

Much as for virtual knots, if a virtual spatial graph diagram has a realizable Gauss code, then that
diagram is equivalent to a classical diagram.  The proofs of Theorems \ref{T:virtualrealize} and
\ref{T:virtualgauss} are similar to those for the analogous results for virtual knots.

\begin{thm} \label{T:virtualgauss}
If two virtual graph diagrams have the same Gauss code, then they
are virtually equivalent.
\end{thm}

\begin{cor}
If a virtual graph diagram has a realizable Gauss code, then it is virtually equivalent to a classical
graph diagram.
\end{cor}

The inclusion of knot theory into virtual knot theory is known to be
injective.  It is not yet known if the same is true for the
inclusion of spatial graph theory into virtual spatial graph theory.

\begin{conj}
If two classical spatial graphs are virtually equivalent, then they
are classically equivalent.
\end{conj}

\section{Intrinsically Virtually Linked Graphs}
\label{S:ivl}

A graph G is called \emph{intrinsically virtually linked of degree
n} ($n$-IVL) if every virtual diagram of G with at most $n$ virtual
crossings contains a non-trivial virtual link whose components are
disjoint cycles in G. Let $IVL_n$ denote the set of graphs that are
intrinsically virtually linked of degree $n$.  We may define
\emph{intrinsically virtually knotted of degree $n$} and $IVK_{n}$
in the same way. These definitions give rise to the natural
filtrations below. Note that $IVL_0$ ($IVK_0$) is simply the set of
classically intrinsically linked (knotted) graphs, so the filtration
provides information about this classical problem.

$$IVL_0 = IVL_{1} \supsetneq IVL_{2} \supsetneq IVL_{3} \supsetneq IVL_{4} \supsetneq
\ldots$$

$$IVK_{0} \supset IVK_{1} \supsetneq IVK_{2} \supset IVK_{3} \supsetneq IVK_{4} \supset
\ldots$$

It is relatively easy to produce examples that show the filtrations are decreasing in the manner
shown (see \cite{fm3}), but the equality $IVL_0 = IVL_1$ is somewhat surprising.

\begin{thm} \label{T:0IVL=1IVL}
If G is intrinsically linked, then G is also intrinsically virtually
linked of degree 1, so $IVL_0 = IVL_1$.
\end{thm}

The proof of Theorem \ref{T:0IVL=1IVL} relies on the fact that every embedding of every
intrinsically linked graph contains a non-split link with odd linking number, and that the
half-integer linking number is an invariant for virtual links.

\begin{figure} [ht]
    $$\includegraphics{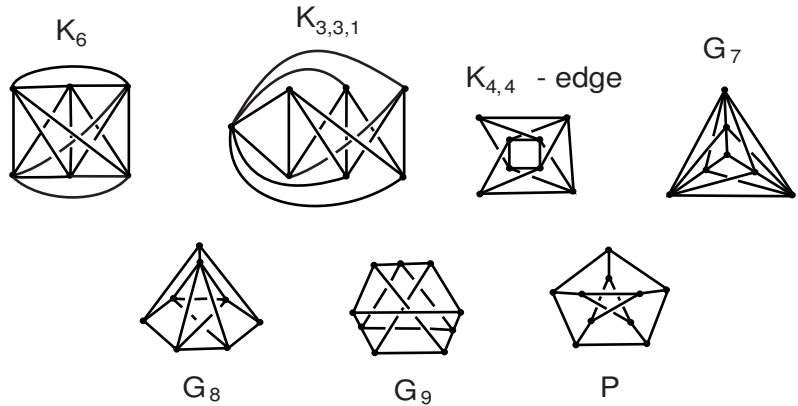}$$
    \caption{The Petersen family of graphs.} \label{F:petersenfamily}
    \end{figure}

\begin{figure} [ht]
    $$\includegraphics{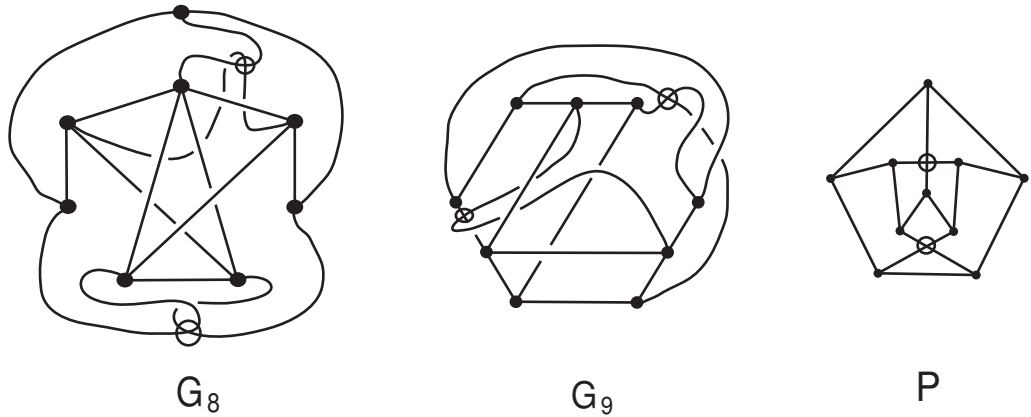}$$
    \caption{Linkless virtual diagrams of $G_{8}$, $G_{9}$ and $P$.} \label{F:2IVL}
    \end{figure}

Every intrinsically linked graph contains one of the Petersen graphs as a minor \cite{r&s}. That
is, the Petersen family, shown in
Figure \ref{F:petersenfamily}, is the minor minimal set for intrinsic linking.  We have classified
the Petersen family of graphs with respect to
the filtration, and clearly all are 1-IVL. Embeddings of $G_{8}$, $G_{9}$ and $P$ with two
virtual crossings that contain no nontrivial links are shown in Figure \ref{F:2IVL}, and hence are
not 2-IVL. The graphs $K_{6}$, $K_{3,3,1}$, $K_{4,4} \setminus e$ and $G_{7}$ are 2-IVL,
but not 3-IVL. These graphs are minor minimal for intrinsic linking, so they are minor minimal
with respect to intrinsic virtual linking of degree 2.  However, the full set of minor minimal
graphs for intrinsic virtual linking of degree 2 is not yet known.

If a knot has non-trivial Jones polynomial, then converting a single classical crossing in any
projection of that knot to a virtual crossing gives a non-trivial virtual knot \cite{fm3}.  Since
every embedding of every known intrinsically knotted graph contains a knot with non-trivial
Jones polynomial, all known 0-IVK graphs are 1-IVK.  Indeed, we conjecture that $IVK_{0} =
IVK_{1}$.  To further study this conjecture, as well as understand intrinsic virtual knotting, we
introduce the notion of \emph{virtual unknotting number}.

Given a diagram $D$ of a virtual knot $K$, define $vu_D(K)$ to be the minimum number of
classical crossings in $D$ which need to be virtualized in order to unknot K.  The virtual
unknotting number of $K$ $vu(K)$ is the minimum of $\{vu_D(K)\}$, taken over all diagrams
D for K.  We have been able to determine the virtual unknotting number for certain knots.

\begin{thm} \label{T:twist}
The virtual unknotting number of any twist knot is 2.
\end{thm}

The virtual unknotting number may be related in some way to the classical unknotting number,
but this is not yet understood.  In this area, as with all of virtual spatial graph theory, many open
problems remain.

\small

\normalsize


\begin{thebibliography}{99}
\bibitem{c&g} J. H. Conway and C. McA. Gordon, Knots and links in spatial graphs, {\it J.
Graph Th.} v. 7, 1983 446-453
\bibitem{fo}  de Fraysseix, H. and Ossona de Mendez, P.:  On a Characterization of Gauss
Codes, {\it
Discrete Comput. Geom.}, v. 22, 1999, pp. 287-295
\bibitem{de}  Dehn, M.:  \"{U}ber Kombinatorische Topologie, {\it Acta Math.} 67, 1936, pp.
123-168
\bibitem{fl}  Flapan, E.:  Symmetries of knotted hypothetical molecular graphs, Applications of
graphs in chemistry and physics. \emph{Discrete Appl. Math.} v. 19
1988, no. 1-3, pp. 157-166
\bibitem{fm}  Fleming, T. and Mellor, B.:  Chord Diagrams and Gauss Codes for Graphs,
preprint, 2005
\bibitem{fm2} Fleming, T. and Mellor, B.:  Virtual Spatial Graphs,
preprint, 2005
\bibitem{fm3} Fleming, T. and Mellor, B.:  Intrinsic Linking and Knotting in Virtual
Spatial Graphs, preprint, 2006
\bibitem{fm4} Fleming, T. and Mellor, B.:  Counting Links in
Complete Graphs, preprint, 2006
\bibitem{joyce} Joyce, D.: A Classifying Invariant of Knots, the Knot Quandle, {\it J. Pure Appl.
Algebra} v. 23, 1982 pp. 37--65
\bibitem{ka}  Kauffman, L.:  Virtual Knot Theory, {\it Europ. J. Combinatorics}, v. 20, 1999,
pp. 663-691
\bibitem{ka2}  Kauffman, L.:  Invariants of Graphs in Three-Space, {\it Trans. Amer. Math.
Soc.}, v. 311, no. 2, 1989, pp. 697-710
\bibitem{ma} Manturov, V.: On Invariants of Virtual Links, {\it Acta Appl. Math.}, v.  72,
2002, pp.  295-309
\bibitem{ne} Nelson, S.:  Unknotting virtual knots with Gauss diagram forbidden moves, {\it J.
Knot Theory Ramif.}, v. 10, no. 6, 2001, pp. 931-935
\bibitem{rr}  Read, R.C. and Rosenstiehl, P.:  On the Gauss Crossing Problem, {\it Colloq. Math.
Soc. Janos
Bolyai}, v. 18, {\it Combinatorics}, Keszthely, Hungary, 1976, pp.
843-876
\bibitem{r&s} N. Robertson, P. Seymour, R. Thomas, Sachs' Linkless Embedding
Conjecture, {\it J. Comb Theory Ser. B} v. 64, 1995, 185-277
\bibitem{sachs} H. Sachs, On Spatial Representations of Finite Graphs, in: A. Hajnal, L.
Lovasz, V.T. S\'os (Eds.), {\it Colloq. Math. Soc. J\'anos Bolyai},
Vol. 37, North-Holland, Amsterdam, 1984, 649-662
\bibitem{dw}  Sumners, D. W.:  Lifting the curtain: using topology to probe the hidden action of
enzymes, {\it Notices Amer. Math. Soc.} v. 42, no. 5, 1995 pp.
528--537
\bibitem{ya}  Yamada, S.:  An Invariant of Spatial Graphs, {\it J. Graph Theory}, v. 13, no. 5,
1989, pp. 537-551
\end{thebibliography}
\end{document}